\documentclass[11pt, reqno, psamsfonts]{amsart}
\pdfoutput=1

\usepackage{amssymb}
\usepackage{amsthm}
\usepackage{amsmath}
\usepackage{latexsym}
\usepackage[T1]{fontenc}
\usepackage[utf8]{inputenc}
\usepackage[russian, french, english]{babel}

\usepackage{graphicx}
\usepackage{wrapfig}
\usepackage{mathtools}
\usepackage{amsbsy}
\usepackage[inline]{enumitem}
\usepackage{mathrsfs}
\usepackage{array}
\usepackage{multicol}
\usepackage{stmaryrd}
\usepackage{cancel}
\usepackage{lmodern}
\usepackage[matha,mathx]{mathabx}
\usepackage{upgreek}
\usepackage{titlesec}
\usepackage{titletoc}
\usepackage[spacing=true,kerning=true,babel=true,tracking=true]{microtype}
\usepackage[shortcuts]{extdash}
\usepackage[foot]{amsaddr}
\usepackage[left=1in,right=1in,top=1in,bottom=1in,bindingoffset=0cm]{geometry}
\usepackage{bm}
\usepackage{centernot}
\usepackage{mdframed}
\usepackage[hidelinks]{hyperref}
\usepackage{xspace}
\usepackage[skins]{tcolorbox}
\usepackage{framed}
\usepackage[justification=centering, labelfont=bf]{caption}
\usepackage{tikz}
\usetikzlibrary{shapes,snakes}
\usetikzlibrary{arrows.meta}
\usetikzlibrary{decorations.pathmorphing}
\usetikzlibrary{patterns}
\usetikzlibrary{shapes.multipart}
\usepackage{float}

\usepackage[
    sortcites,
    backend=biber, style=alphabetic, sorting=nyt, maxnames=100,backref=true]{biblatex}
\title{\sffamily {On strongly and robustly critical graphs}}
\date{}

\author{Anton~Bernshteyn}
\address{\normalfont (AB) Department of Mathematics, University of California, Los Angeles, CA, USA}
\email{bernshteyn@math.ucla.edu}

\author{Hemanshu~Kaul}
\address{\normalfont (HK) Department of Applied Mathematics, Illinois Institute of Technology, Chicago, IL, USA}
\email{kaul@iit.edu}

\author{Jeffrey~A.~Mudrock}
\address{\normalfont (JM) Department of Mathematics and Statistics, University of South Alabama, Mobile, AL, USA}
\email{mudrock@southalabama.edu}

\author{Gunjan~Sharma}
\address{\normalfont (GS) Department of Applied Mathematics, Illinois Institute of Technology, Chicago, IL, USA}
\email{gsharma7@hawk.iit.edu}

\thanks{AB's research is partially supported by the NSF CAREER grant DMS-2528522 and the Sloan Research Fellowship.}

\newtheoremstyle{bfnote}%
{}{}%
{\slshape}{}%
{\bfseries}{\bfseries.}%
{ }%
{\thmname{#1}\thmnumber{ #2}\thmnote{ \ep{\normalfont{}#3}}}

\theoremstyle{bfnote}
\newtheorem{thm}{Theorem}[section]
\newtheorem*{thm*}{Theorem}
\newtheorem{pro}[thm]{Proposition}
\newtheorem{lem}[thm]{Lemma}
\newtheorem{cl}[thm]{Claim}

\newtheorem{cor}[thm]{Corollary}
\newtheorem{conj}[thm]{Conjecture}

\newtheorem*{cor*}{Corollary}

\theoremstyle{definition}
\newtheorem{defn}[thm]{Definition}
\newtheorem*{defn*}{Definition}
\newtheorem{exmp}[thm]{Example}

\newtheorem*{exmp*}{Example}
\newtheorem{prob}[thm]{Problem}

\theoremstyle{remark}
\newtheorem*{ques*}{Question}
\newtheorem*{remk*}{Remark}

\newcommand*{\myproofname}{Proof}
\newenvironment{claimproof}[1][\myproofname]{\begin{proof}[#1]}{\end{proof}}

\makeatletter
\newcommand{\neutralize}[1]{\expandafter\let\csname c@#1\endcsname\count@}
\makeatother

\newenvironment{theocopy}[1]
{\renewcommand{\thethm}{\ref{#1}}%
	\neutralize{thm}\phantomsection
	\begin{thm}}
	{\end{thm}}

\newenvironment{propcopy}[1]
{\renewcommand{\thethm}{\ref{#1}}%
	\neutralize{thm}\phantomsection
	\begin{pro}}
	{\end{pro}}

\newcommand{\0}{\varnothing}
\newcommand{\set}[1]{\{#1\}}
\newcommand{\N}{{\mathbb{N}}}

\renewcommand{\epsilon}{\varepsilon}

\renewcommand{\phi}{\varphi}
\renewcommand{\theta}{\vartheta}
\renewcommand{\leq}{\leqslant}
\renewcommand{\geq}{\geqslant}
\newcommand{\defeq}{\coloneqq}

\newcommand{\bemph}[1]{{\normalfont#1}} 
\newcommand{\ep}[1]{\bemph{(}#1\bemph{)}} 

\newcommand{\emphd}[1]{{\fontseries{b}\selectfont\textsf{#1}}}
\newcommand{\dom}{\mathsf{dom}}

\newcommand{\rest}[2]{{{#1}\vert_{#2}}}

\renewcommand{\deg}{\mathsf{deg}}

\numberwithin{equation}{section}

\newenvironment{scproof}[1][]{\begin{proof}[\textsc{\upshape{Proof}}#1]}{\end{proof}}

\titleformat{\section}[block]{\large\bfseries\sffamily}{\thesection.}{1ex}{}
\titleformat{\subsection}[block]{\bfseries\sffamily}{\thesubsection.}{1ex}{}
\titleformat{\subsubsection}[block]{\itshape}{\bfseries\upshape\sffamily\thesubsubsection.}{1ex}{}

\titlespacing*{\section}{0pt}{*3}{*1}
\titlespacing*{\subsection}{0pt}{*3}{*1}
\titlespacing*{\subsubsection}{0pt}{*2}{*1}

\titlecontents{section}
 [1.5em] %
 {\smallskip}
 {\bfseries\thecontentslabel\hspace{1.02em}}
{\bfseries}
 {\,\,\titlerule*[0.77pc]{}\bfseries\contentspage}
\titlecontents{subsection}
 [4em] %
 {\smallskip}
 {\thecontentslabel\hspace{1.02em}}
{\hspace*{2.32em}}
 {\,\,\titlerule*[0.77pc]{.}\contentspage}

\renewbibmacro{in:}{}

\renewbibmacro*{volume+number+eid}{%
	\printfield{volume}%
	\setunit*{\addnbspace}
	\printfield{number}%
	\setunit{\addcomma\space}%
	\printfield{eid}}

\DeclareFieldFormat[article]{volume}{\textbf{#1}\space}
\DeclareFieldFormat[article]{number}{\mkbibparens{#1}}

\DeclareFieldFormat{journaltitle}{#1,}
\DeclareFieldFormat[thesis]{title}{\mkbibemph{#1}\addperiod}
\DeclareFieldFormat[article, unpublished, thesis]{title}{\mkbibemph{#1},}
\DeclareFieldFormat[book]{title}{\mkbibemph{#1}\addperiod}
\DeclareFieldFormat[unpublished]{howpublished}{#1, }

\DeclareFieldFormat{pages}{#1}

\DeclareFieldFormat[article]{series}{Ser.~#1\addcomma}

\setlist{topsep=3pt,itemsep=3pt}

\newenvironment{graph}[1][scale=1]{
\begin{tikzpicture}[#1]
\tikzstyle{vertex}=[circle, draw, fill, inner sep=0pt, minimum size=4pt]%
\tikzstyle{every path}=[line width=0.5pt]%
\tikzstyle{G}=[dashed]%
\tikzstyle{F}=[solid]
}{\end{tikzpicture}}

\begin{document}

\maketitle



\begin{abstract}
    In extremal combinatorics, it is common to focus on structures that are minimal with respect to a certain property. In particular, critical and list-critical graphs occupy a prominent place in graph coloring theory. Stiebitz, Tuza, and Voigt introduced strongly critical graphs, i.e., graphs that are $k$-critical yet $L$-colorable with respect to every non-constant assignment $L$ of lists of size $k-1$. Here we strengthen this notion and extend it to the framework of DP-coloring (or correspondence coloring) by defining robustly $k$-critical graphs as those that are not $(k-1)$-DP-colorable, but only due to the fact that $\chi(G) = k$. We then seek general methods for constructing robustly critical graphs. Our main result is that if $G$ is a critical graph (with respect to ordinary coloring), then the join of $G$ with a sufficiently large clique is robustly critical; this is new even for strong criticality. 


\medskip


\smallskip


\end{abstract}

\section{Introduction}\label{intro}

    All graphs in this paper are finite and simple. A useful approach in graph theory is to focus one's attention on graphs that are minimal with respect to a property of interest. In the study of graph coloring, this philosophy was applied by Dirac \cite{Dirac1,Dirac2}, who introduced the notion of a critical graph. A graph $G$ is called \emphd{$k$-critical} if $\chi(G) = k$ and $\chi(G') < k$ for every proper subgraph $G'$ of $G$; a graph $G$ is \emphd{critical} if it is $k$-critical for some $k$. Every graph with chromatic number $k$ has a $k$-critical subgraph; therefore, insight into the structure of critical graphs has the potential to shed light on the coloring properties of graphs in general. For example, the structure of critical graphs can be used to derive strong conclusions concerning chromatic numbers of graphs embedded in a given surface; see, e.g., \cite{Dirac3,KY,LiuPostle,MSR,Dirac5,Gallai,Thomassen1,Thomassen2}. 

    Another fundamental concept in graph theory, \emphd{list-coloring}, was introduced in the 1970s by Vizing \cite{Viz} and, independently, Erd\H{o}s, Rubin, and Taylor \cite{ERT}; for textbook introductions, see 
    \cites[\S14.5]{BondyMurty}[\S5.4]{Diestel}[\S8.4]{W01}. In the list-coloring  framework, each vertex $v$ of a graph $G$ is assigned a set $L(v)$, called its \emphd{list of available colors}. The objective is to choose for each $v \in V(G)$ a color $f(v) \in L(v)$ so that adjacent vertices receive different colors; such a coloring $f$ is called a \emphd{proper $L$-coloring} of $G$. Note that if $L(v) = [k]$ 
    for every vertex $v \in V(G)$, this turns into the ordinary $k$-coloring problem. If a proper $L$-coloring of $G$ exists, we say $G$ is \emphd{$L$-colorable}; otherwise, we say that $L$ is a \emphd{bad} list assignment for $G$. If $|L(v)| = k$ for all $v \in V(G)$, we call $L$ a \emphd{$k$-assignment}. The \emphd{list-chromatic number} (also called the \emphd{choosability}) of $G$, denoted by $\chi_\ell(G)$, is the minimum $k$ such that $G$ is $L$-colorable for every $k$-assignment $L$. It is immediate from the definition that $\chi_\ell(G) \geq \chi(G)$. This inequality can be strict, as, for example, $\chi_\ell(K_{n,n}) = \Theta(\log n)$ while $\chi(K_{n,n}) = 2$ \cite{ERT}.
 
    Following the ``minimal example'' philosophy, it is natural to consider critical graphs with respect to list-coloring, and indeed they have been researched extensively with impressive results paralleling those in the study of critical graphs for ordinary coloring; for example, see \cite{BMS,PostleThomas,DeVKM} for some applications of list-criticality in the theory of embedded graphs. Specifically, a graph $G$ is called \emphd{$L$-critical}, where $L$ is a list assignment for $G$, if $L$ is bad for $G$ but every proper subgraph of $G$ is $L$-colorable.\footnote{Here we employ a standard abuse of terminology and say that, given a list assignment $L$ for $G$, a subgraph $G' \subseteq G$ is $L$-colorable if it is $L'$-colorable, where $L'$ is the list assignment for $G'$ defined by $L'(v) \defeq L(v)$ for all $v \in V(G')$.} We emphasize that this definition explicitly depends on the list assignment $L$; for example, an $L$-critical graph $G$ with respect to a $k$-assignment $L$ may have a proper subgraph $G'$ that is $L'$-critical for some other $k$-assignment $L'$ \cite[Example 3]{STV09}.

    By definition, a graph $G$ is $k$-critical if and only if it is $L$-critical for the list assignment $L$ given by $L(v) = [k-1]$ for all $v \in V(G)$. This means that if we remove any vertex or edge from $G$, the resulting graph becomes $L$-colorable. What if we modify the list assignment $L$ itself? If the constant list assignment $L$ 
    is the only (up to renaming the colors) bad $(k-1)$-assignment for $G$, then $G$ is called strongly critical; this notion was introduced by Stiebitz, Tuza, and Voigt \cite{STV09} in the course of their general investigation into the structure of list-critical graphs.

    \begin{defn}[Strongly critical graphs]\label{defn:strong}
        A graph $G$ is \emphd{strongly $k$-critical}\footnote{In \cite{STV09}, the term ``strong $k$-critical'' was used.} for $k \geq 1$ if $G$ is $k$-critical and every bad $(k-1)$-assignment for $G$ is constant (i.e., it assigns the same $(k-1)$-element set of colors to every vertex). If $G$ is strongly $k$-critical for some $k$, we say that $G$ is \emphd{strongly critical}.
    \end{defn}

    In other words, for a  strongly $k$-critical graph $G$, its chromatic number is the only obstruction to $(k-1)$-list-coloring. The examples described below are taken from \cite[\S2.2]{STV09}.

    \begin{figure}[t]
			\centering
			\begin{tikzpicture}[scale=0.85, every node/.style={transform shape}]
                \node[fill=black!15!white,ellipse,minimum width=0.5cm,minimum height=1.5cm] at (0,0) {};
                \node[fill=black!15!white,ellipse,minimum width=0.5cm,minimum height=2.5cm] at (1.5,0) {};
                \node[fill=black!15!white,ellipse,minimum width=0.5cm,minimum height=1.5cm] at (3,0) {};
                \node[fill=black!15!white,ellipse,minimum width=0.5cm,minimum height=2.5cm] at (4.5,0) {};
   
			    \filldraw (0,-0.5) circle (2pt);
                \filldraw (0,0.5) circle (2pt);

                \filldraw (1.5,1) circle (2pt);
                \filldraw (1.5,0) circle (2pt);
                \filldraw (1.5,-1) circle (2pt);

                \filldraw (3,-0.5) circle (2pt);
                \filldraw (3,0.5) circle (2pt);

                \filldraw (4.5,1) circle (2pt);
                \filldraw (4.5,0) circle (2pt);
                \filldraw (4.5,-1) circle (2pt);

                \draw (0,-0.5) -- (0,0.5);
                \draw (1.5,1) -- (1.5,0) -- (1.5,-1) to[bend left=45] (1.5,1);

                \draw (3,-0.5) -- (3,0.5);
                \draw (4.5,1) -- (4.5,0) -- (4.5,-1) to[bend left=45] (4.5,1);

                \draw (0,-0.5) -- (1.5,1) -- (0,0.5) -- (1.5,0) -- (0,-0.5) -- (1.5,-1) -- (0,0.5);

                \draw (3,-0.5) -- (1.5,1) -- (3,0.5) -- (1.5,0) -- (3,-0.5) -- (1.5,-1) -- (3,0.5);

                \draw (3,-0.5) -- (4.5,1) -- (3,0.5) -- (4.5,0) -- (3,-0.5) -- (4.5,-1) -- (3,0.5);

                \filldraw (2.25,-2) circle (2pt);

                \draw (2.25,-2) to[out=180,in=-90] (0,-0.5);
                \draw (2.25,-2) to[out=180,in=200,looseness=2] (0,0.5);
                \draw (2.25,-2) to[out=0,in=-120] (4.5,-1);
                \draw (2.25,-2) to[out=0,in=-20,looseness=2] (4.5,0);
                \draw (2.25,-2) to[out=0,in=-30,looseness=2.5] (4.5,1);

                \node at (2.25,-2.3) {$z$};
                \node at (0,1) {$X_1$};
                \node at (1.5,1.5) {$X_2$};
                \node at (3,1) {$Y_2$};
                \node at (4.5,1.5) {$Y_1$};
                
			\end{tikzpicture}
   \caption{The strongly (and robustly) $6$-critical graph $E_{6,2,3}$.}\label{fig:E623}
	\end{figure}
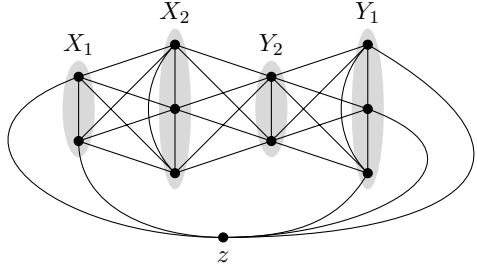

    \begin{exmp}[The graphs $E_{k,a,b}$]\label{exmp:E}
        It is easy see that complete graphs and odd cycles are strongly critical. Here is a more interesting family of examples. For integers $k \geq 3$ and $1 \leq a$, $b \leq k - 2$ such that $a+b \geq k-1$, we define a graph $E_{k,a,b}$ as follows. Let $X_1$, $X_2$, $Y_1$, $Y_2$ be disjoint sets with $|X_1| = a$, $|Y_1| = b$, $|X_2| = k - 1 - a$, and $|Y_2| = k - 1 - b$, and let $z$ be an additional vertex. The vertex set of $E_{k,a,b}$ is $V(E_{k,a,b}) \defeq X_1 \cup X_2 \cup Y_1 \cup Y_2 \cup \set{z}$. The sets $X \defeq X_1 \cup X_2$ and $Y \defeq Y_1 \cup Y_2$ are cliques in $E_{k,a,b}$, the neighborhood of $z$ is $X_1 \cup Y_1$, and a vertex $x \in X$ is adjacent to $y \in Y$ if and only if $x \in X_2$ and $y \in Y_2$ (see Figure~\ref{fig:E623} for an illustration). Then the graph $E_{k,a,b}$ is strongly $k$-critical \cite[Example 8]{STV09}. In the special case $b = k-2$, $E_{k,a,k-2}$ is called a \emphd{Dirac graph} \cite[Example 7]{STV09}. Dirac graphs play an important role in the theory of critical and list-critical graphs due to their relatively low number of edges \cite{Dirac4,BK18}.
    \end{exmp}

    Stiebitz, Tuza, and Voigt also observed that the join of a strongly critical graph and a clique is itself strongly critical; here the \emphd{join} of graphs $G$ and $H$, denoted by $G \vee H$, is the graph obtained from vertex-disjoint copies of $G$ and $H$ by making every vertex of $G$ adjacent to every vertex of $H$. 

    \begin{pro}[{\cite[Proposition 9]{STV09}}]\label{prop:join}
        If $G$ is a strongly critical graph, then $G \vee K_t$ is also strongly critical for any $t \in \N$.
    \end{pro}

    To the best of our knowledge, the only examples of strongly critical graphs that have appeared in the literature are complete graphs, $C_{2\ell + 1} \vee K_t$, and $E_{k,a,b} \vee K_t$, obtained by combining Example~\ref{exmp:E} with Proposition~\ref{prop:join}. This should be contrasted with the richness of the class of critical graphs \cite{rich1,rich2,rich3,rich4,rich5,rich6,rich7,rich8,rich9,GL74}; see also the survey \cite{rich_survey} by Sachs and Stiebitz. For example, Greenwell and Lov\'{a}sz \cite[Theorem 3]{GL74} proved that every graph $G$ is an induced subgraph of some critical graph.
    
    Here we contribute a new general way of building strongly critical graphs. Namely, we show that for any critical graph $G$, its join with a sufficiently large clique becomes strongly critical:
    
    \begin{thm}[Joins with complete graphs make critical graphs strongly critical]\label{theo:joins_make_strong}
        If $G$ is a critical graph with $m$ edges, then for all $t \geq 3m$, the graph $G \vee K_t$ is strongly critical.
    \end{thm}

    The assumption that $G$ is critical in Theorem~\ref{theo:joins_make_strong} is necessary, since otherwise $G \vee K_t$ would fail to be critical (let alone strongly critical), no matter how large $t$ is. Note that every strongly critical graph $G$ is \emphd{chromatic-choosable}, i.e., $\chi_\ell(G) = \chi(G)$. Therefore, Theorem~\ref{theo:joins_make_strong} is a strengthening for critical graphs $G$ of a result of Ohba that for all large $t$, $G \vee K_t$ is chromatic-choosable \cite{Ohba}.


 
    A graph $G$ is \emphd{$k$-vertex-critical} if $\chi(G) = k$ and $\chi(G') < k$ for every proper \emph{induced} subgraph $G'$ of $G$; if $G$ is $k$-vertex-critical for some $k$, we say that $G$ is \emphd{vertex-critical}. A variant of Definition~\ref{defn:strong} for vertex-criticality was given by the second and third named authors in \cite{KM18}:

    \begin{defn}[Strongly chromatic-choosable graphs]
        A graph $G$ is \emphd{strongly $k$\-/chromatic\-/choosable} for $k \geq 1$ if it is $k$-vertex-critical and every bad $(k-1)$-assignment for $G$ is constant. If $G$ is strongly $k$-chromatic-choosable for some $k$, we say that $G$ is \emphd{strongly chromatic-choosable}.
    \end{defn}

    We prove a version of Theorem~\ref{theo:joins_make_strong} for vertex-critical graphs:

    \begin{thm}[Vertex version of Theorem~\ref{theo:joins_make_strong}]\label{theo:vertex}
        If $G$ is a vertex-critical graph with $m$ edges, then for all $t \geq 3m$, the graph $G \vee K_t$ is strongly chromatic-choosable.
    \end{thm}



    By definition, every strongly critical graph is strongly chromatic-choosable. On the other hand, in \cite[\S2.1]{KM18} it is shown that there exist strongly chromatic-choosable graphs that fail to be strongly critical. Theorem~\ref{theo:vertex} provides further such examples: take any graph $G$ that is vertex-critical but not critical (see \cite{vtx-critical} for a recent paper by Martinsson and Steiner that gives constructions of vertex-critical graphs that are far from being critical) and consider $G \vee K_t$ for sufficiently large $t$.
    
    Note that Theorem~\ref{theo:vertex} implies Theorem~\ref{theo:joins_make_strong}. Indeed, if $G$ is a critical graph with $m$ edges, then, by Theorem~\ref{theo:vertex}, $G \vee K_t$ is strongly chromatic-choosable for any $t \geq 3m$. Moreover, it is easy to see that the join of a critical graph and a clique is critical. Hence, $G \vee K_t$ is strongly critical, as desired.

    A further generalization of list-coloring that has been a topic of active research in recent years is so-called \emphd{DP-coloring} (also known as \emphd{correspondence coloring}), which was invented by Dvo\v{r}\'ak and Postle \cite{DP15} and is closely related to \emph{local conflict coloring} introduced by Fraigniaud, Heinrich, and Kosowski \cite{FHK}. Even though DP-coloring has only emerged relatively recently, it has already garnered considerable attention.\footnote{According to MathSciNet, the paper \cite{DP15} by Dvo\v{r}\'ak and Postle has over 100 citations at the time of writing.} In the DP-coloring framework, not only the lists of available colors but also the identifications between them may vary from edge to edge. The way the correspondences between the colors are arranged can be conveniently captured by an additional structure called a cover of $G$:

    \begin{defn}[Covers and DP-colorings]
        A \emphd{cover} of a graph $G$ is a pair $\mathcal{H} = (L,H)$, where:
        \begin{itemize}
            \item $H$ is a graph and $L$ is a function assigning to each $v \in V(G)$ a subset $L(v) \subseteq V(H)$,

            \item the sets $L(v)$ for $v \in V(G)$ are disjoint, independent\footnote{Some sources, e.g., \cite{BKP}, require the sets $L(v)$ to be cliques instead. This distinction makes no difference for the way the definition is used.} in $H$, and satisfy $V(H) = \bigcup_{v \in V(G)} L(v)$, 

            \item if $E_H(L(u),L(v)) \neq \0$, then $uv \in E(G)$ and $E_H(L(u), L(v))$ is a matching.
        \end{itemize}
        We stress that the matchings between $L(u)$ and $L(v)$ for $uv \in E(G)$ need not be perfect (and may even be empty). 
        The vertices of $H$ are referred to as \emphd{colors}. An \emphd{independent transversal} of $\mathcal{H}$ is an independent set $T \subseteq V(H)$ in $H$ containing exactly one vertex from each list $L(v)$. If $\mathcal{H}$ has an independent transversal $T$, we call $T$ a \emphd{proper $\mathcal{H}$-coloring} of $G$ and say that $G$ is \emphd{$\mathcal{H}$-colorable}. If $G$ is not $\mathcal{H}$-colorable, we call $\mathcal{H}$ a \emphd{bad} cover of $G$.

        A cover $\mathcal{H} = (L,H)$ of $G$ is called \emphd{$k$-fold} if $|L(v)| = k$ for all $v \in V(G)$. (Here $k$ is allowed to be $0$. In the unique $0$-fold cover $(L,H)$ of $G$, $H$ is the empty graph and $L(v) = \0$ for all $v \in V(G)$.) The \emphd{DP-chromatic number} of $G$, denoted by $\chi_\mathsf{DP}(G)$, is the minimum $k$ such that $G$ is $\mathcal{H}$-colorable for every $k$-fold cover $\mathcal{H}$.
    \end{defn}

    The ordinary $k$-coloring problem for a graph $G$ is equivalent to $\mathcal{H}$-coloring with respect to the cover $\mathcal{H} = (L, H)$ where $L(v) = \set{v_1, \ldots, v_k}$ for each $v \in V(G)$ and $E(H) \defeq \set{u_iv_i \,:\, uv \in E(G), \, i \in [k]}$. More generally, from any $k$-assignment $L$ for $G$, one can construct an associated $k$-fold cover $\mathcal{H}_L$ and a one-to-one correspondence between proper $L$-colorings and proper $\mathcal{H}_L$-colorings of $G$ (in particular, $G$ is $L$-colorable if and only if it is $\mathcal{H}_L$-colorable) \cite{DP15}. It follows that $\chi_{\mathsf{DP}}(G) \geq \chi_\ell (G)$ for every graph $G$. This inequality can be strict, as exemplified by the fact that $\chi_\ell(K_{n,n}) = \Theta(\log n)$ but $\chi_\mathsf{DP}(K_{n,n}) = \Theta(n/\log n)$ \cite{B}.

    By analogy with Definition~\ref{defn:strong}, we say that a graph $G$ is {robustly $k$-critical} if it is $k$-critical and the only bad $(k-1)$-cover of $G$ is the one that encodes the usual $(k-1)$-coloring problem; such covers are called {canonical}. Here is the formal definition:

    \begin{defn}[Canonical covers and robustly critical graphs]
        A $k$-fold cover $\mathcal{H} = (L,H)$ of a graph $G$ is \emphd{canonical} if it admits a \emphd{canonical labeling}, i.e., a mapping $\lambda \colon V(H) \to [k]$ such that
        \begin{itemize}
            \item for each $v \in V(G)$, $\rest{\lambda}{L(v)}$ is a bijection from $L(v)$ to $[k]$, and
            \item for all $uv \in E(G)$ and $c \in L(u)$, $c' \in L(v)$, we have $cc' \in E(H)$ if and only if $\lambda(c) = \lambda(c')$. 
        \end{itemize}
        Note that, in particular, the $0$-fold cover of $G$ is canonical.
    
        A graph $G$ is \emphd{robustly $k$-critical} for $k \geq 1$ if it is $k$-critical and every bad $(k-1)$-fold cover of $G$ is canonical. If $G$ is robustly $k$-critical for some $k$, we say that $G$ is \emphd{robustly critical}.
    \end{defn}

    \begin{figure}[t]
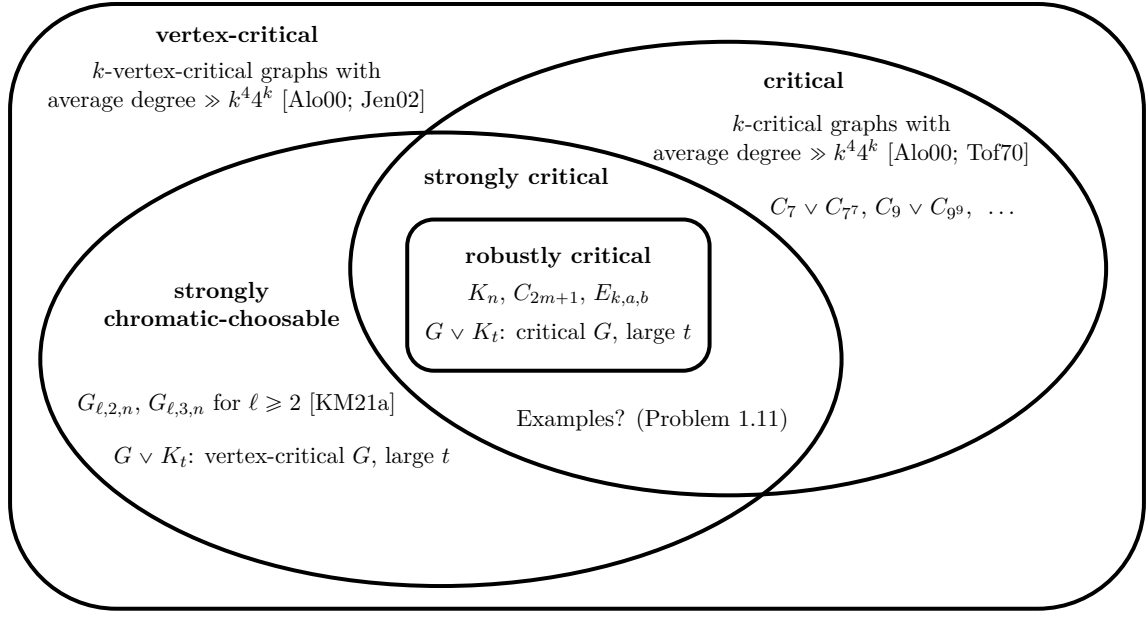

\centering
\begin{graph}

\draw[line width=1.5,rounded corners=40] (-7.5,-4.5) rectangle (7.5,3.5);
\node [scale =0.8] at (-4.5,3.1) {\textbf{vertex-critical}};
\node [scale =0.8] at (-4.5,2.4) {\begin{tabular}{c}$k$-vertex-critical graphs with \\ average degree $\gg k^44^k$ \cite{rich7,Alon}\end{tabular}};

\draw[line width=1.5] (2,0) ellipse (5cm and 3cm);
\node [scale =0.8] at (3,2.5) {\textbf{critical}};
\node [scale =0.8] at (3.5,1.7) {\begin{tabular}{c}$k$-critical graphs with \\ average degree $\gg k^44^k$ \cite{rich3,Alon}\end{tabular}};
\node [scale =0.8] at (4.2,0.8) {$C_{7} \vee C_{7^7}$, $C_9 \vee C_{9^9}$, \ \ldots};


\draw[line width=1.5] (-1.8,-1.2) ellipse (5.3cm and 3cm);


\node [scale =0.8] at (-4.7,-0.5) {\begin{tabular}{c}\textbf{strongly} \\ \textbf{chromatic-choosable}\end{tabular}};
\node [scale =0.8] at (-4.5,-1.8) {$G_{\ell,2,n}$, $G_{\ell, 3, n}$ for $\ell \geq 2$ \cite{KM18}};
\node [scale =0.8] at (-3.9,-2.5) {$G \vee K_t$: vertex-critical $G$, large $t$};

\node [scale =0.8] at (-0.8,1.2) {\textbf{strongly critical}};

\begin{scope}[yshift=-10,xshift=-7]
\draw[line width=1.5,rounded corners=10] (-2,-1) rectangle (2,1);
\node [scale =0.8] at (0,0.5) {\textbf{robustly critical}};
\node [scale =0.8] at (0,0) {$K_{n},\, C_{2m+1},\, E_{k,a,b}$};
\node [scale =0.8] at (0,-0.5) {$G \vee K_t$: critical $G$, large $t$};
\end{scope}


\node [scale =0.8] at (1,-2) {Examples? (Problem~\ref{prob:examples})};

\end{graph}
\caption{A Venn diagram for various notions of criticality.}
\label{fig:venn}
\end{figure}

    Figure~\ref{fig:venn} provides a diagram of the relationship between criticality, vertex-criticality, strong chromatic-choosability, strong criticality, and robust criticality and mentions a few examples of graphs in the different categories. In particular, since DP-coloring generalizes list-coloring, it is clear that every robustly critical graph is strongly critical (but we do not know whether the converse holds---see Problem~\ref{prob:examples} below). It is also easy to check that complete graphs and odd cycles are robustly critical \cite{KMG21}. The graphs in Example~\ref{exmp:E} are robustly critical and Proposition~\ref{prop:join} holds in the context of robust criticality:

    \begin{pro}[Examples of robustly critical graphs]\label{prop:examples}
        \mbox{}
        
        \begin{enumerate}[label=\ep{\normalfont\roman*}]
            \item\label{item:Ekab} The graphs $E_{k,a,b}$ from Example~\ref{exmp:E} are robustly critical.

            \item\label{item:join} If $G$ is robustly critical, then so is $G \vee K_t$ for any $t \in \N$.
        \end{enumerate} 
    \end{pro}

    It is worth pointing out that there is no ``vertex version'' of robust criticality analogous to strong chromatic-choosability. This is because if $G$ is a graph with $\chi(G) = k \geq 2$ such that every bad $(k-1)$-fold cover of $G$ is canonical, then $G$ must be $k$-critical. Otherwise, there is an edge $e \in E(G)$ such that $\chi(G - e) = k$, and we can form a non-canonical bad $(k-1)$-fold cover $\mathcal{H} = (L,H)$ of $G$ by taking $L(v) \defeq \set{v_1, \ldots, v_{k-1}}$ for all $v \in V(G)$ and setting $E(H) \defeq \set{u_iv_i \,:\, uv \in E(G) \setminus \set{e}, \, i \in [k - 1]}$ (i.e., the cover is almost canonical except that the matching corresponding to $e$ is empty).  

    In \cite{BK17}, Kostochka, Zhu, and the first named author proved a version of Ohba's theorem \cite{Ohba} for DP-coloring, i.e., they showed that for every graph $G$, $\chi_\mathsf{DP}(G \vee K_t) = \chi(G \vee K_t)$ if $t$ is large enough. We extend Theorem~\ref{theo:joins_make_strong} to the DP-coloring setting and show that if the starting graph $G$ is critical, then for large enough $t$, $G \vee K_t$ is robustly critical:

    \begin{thm}[Joins with larger complete graphs make critical graphs robustly critical]\label{theo:joins_make_robust}
        If $G$ is a critical graph with $m$ edges, then for all $t \geq 100\,m^3$, the graph $G \vee K_t$ is robustly critical.
    \end{thm}



    Notice that while Theorem~\ref{theo:joins_make_robust} yields a stronger conclusion than Theorem~\ref{theo:joins_make_strong} (robust criticality in place of strong criticality), it also requires a higher lower bound on $t$. We do not know if the higher bound is actually necessary. Indeed, we do not have any examples of graphs that are strongly critical but not robustly critical; however, we suspect such examples must exist.

    \begin{prob}\label{prob:examples}
        Give an example of a graph $G$ that is strongly critical but not robustly critical.
    \end{prob}

    Another natural problem is to reduce the lower bound on $t$ in Theorems~\ref{theo:joins_make_strong}, \ref{theo:vertex}, and \ref{theo:joins_make_robust}. The related question of estimating the minimum $t$ such that $\chi_\ell(G \vee K_t) = \chi(G \vee K_t)$ or $\chi_\mathsf{DP}(G \vee K_t) = \chi(G \vee K_t)$ has been studied, e.g., by Ohba \cite{Ohba}, Enomoto,  Ohba,  Ota, and Junko \cite{EOOS}, Kostochka, Zhu, and the first named author \cite{BK17}, and Zhang and Dong \cite{ZhangDong}.

    Yet another intriguing open question concerns enumerative aspects of graph coloring. Let $P(G,k)$ denote the number of proper $k$-colorings of a graph $G$ using the colors $1$, \ldots, $k$. This is a polynomial function of $k$, called the \emphd{chromatic polynomial} of $G$. Now, if $G$ is robustly $k$-critical, that should intuitively mean that the most difficult instance of the DP-coloring problem for $G$ is just ordinary $k$-coloring. Hence, the following question does not seem entirely far-fetched:

    \begin{prob}\label{prob:counting}
        Suppose $G$ is a robustly $k$-critical graph. Does it follow that for every $k$-fold cover $\mathcal{H}$ of $G$, there exist at least $P(G, k)$ proper $\mathcal{H}$-colorings?
    \end{prob}

    Note that Problem~\ref{prob:counting} has a positive answer when $G$ is a complete graph, an odd cycle \cite{KM20}, or, more generally, the join of an odd cycle and a clique \cite{CL1}. On the other hand, the following special case of Problem~\ref{prob:counting} is open:

    \begin{conj}\label{conj:counting}
        Let $G$ be a $k$-critical graph. Then, for all large enough $t$ and for every $(k+t)$-fold cover $\mathcal{H}$ of $G \vee K_t$, the number of proper $\mathcal{H}$-colorings of $G \vee K_t$ is at least $P(G \vee K_t, k+t)$.
    \end{conj}

    The following variant of Problem~\ref{prob:counting} for list-coloring was proposed in \cite[485]{KM18} and also remains open: If $G$ is a strongly $k$-chromatic-choosable graph, does it follow that $G$ has at least $P(G,k)$ proper $L$-colorings for every $k$-assignment $L$ for $G$? And if $q > k$, does $G$ have at least $P(G,q)$ proper $L$-colorings with respect to every $q$-assignment? 
    Problem~\ref{prob:counting} and Conjecture~\ref{conj:counting} belong to the broader research program in the enumeration of DP-colorings initiated in \cite{KM20}: For which graphs $G$ is $P_{\mathsf{DP}}(G,k)$ equal to $P(G,k)$ for $k$ large enough? Here $P_{\mathsf{DP}}(G,k)$ is the DP-coloring analog of the chromatic polynomial $P(G,k)$; see \cite{KM18,KM20} for further discussion.

    In addition to its intrinsic interest, the concept of robust criticality has the potential to yield new bounds on DP-chromatic numbers. For example, in a recent paper \cite{upcoming}, the second, third, and fourth named authors employ robustly critical graphs to obtain new bounds on DP-chromatic numbers of Cartesian products, generalizing and extending previously known results in the list- and DP-coloring setting.

    We conclude the introduction with an outline of the structure of the paper. In \S\ref{sec:basics} we establish some preliminary results, including Proposition~\ref{prop:examples}. In \S\ref{sec:vertex} we verify Theorem~\ref{theo:vertex} (and hence also its corollary, Theorem~\ref{theo:joins_make_strong}). Finally, in \S\ref{sec:robust} we prove Theorem~\ref{theo:joins_make_robust}.

    \section{Some basic facts on robustly critical graphs}\label{sec:basics}

    In this section we establish some basic properties of robustly critical graphs as well as tools that are helpful in proving a graph is robustly critical. We begin by observing that a cover of a robustly $k$-critical graph with lists of size at least $k-1$ can only be bad if all the lists have size exactly $k-1$ (in which case the cover must be canonical).

    \begin{lem}\label{lemma:excess}
        If $G$ is a robustly $k$-critical graph and $\mathcal{H} = (L,H)$ is a bad cover of $G$ such that $|L(u)| \geq k - 1$ for all $u \in V(G)$, then $\mathcal{H}$ is a canonical $(k-1)$-fold cover.
    \end{lem}
    \begin{scproof}
        The statement is clear for $k = 1$, so we may assume $k \geq 2$. In particular, $G$ has no isolated vertices. If $\mathcal{H}$ is a $(k-1)$-fold cover, then we are done by the definition of robust criticality. Now suppose that there is a vertex $u \in V(G)$ with $|L(u)| \geq k$. We will show that $G$ is $\mathcal{H}$-colorable. By removing some of the colors from $H$ if necessary, we may arrange that $|L(u)| = k$ and $|L(v)| = k - 1$ for all $v \neq u$. Take any neighbor $v$ of $u$. Since $|L(v)| = k - 1 < k =  |L(u)|$, there is a color $c \in L(u)$ that has no neighbor in $L(v)$. Form a $(k-1)$-fold cover $\mathcal{H}' = (L', H')$ of $G$ by removing from $L(u)$ any one color other than $c$. The matching corresponding to the edge $uv$ in $\mathcal{H}'$ is not perfect (because $c$ has no neighbor in it), hence $\mathcal{H}'$ is not canonical. This implies that $G$ is $\mathcal{H}'$-colorable, and thus it is also $\mathcal{H}$-colorable, a contradiction.
    \end{scproof}

    A $k$-fold cover $\mathcal{H} = (L,H)$ of a graph $G$ is \emphd{full} if for all $uv \in E(G)$, the matching $E_H(L(u), L(v))$ is perfect. Note that if $\mathcal{H} = (L,H)$ is any $k$-fold cover of $G$, then $\mathcal{H}$ has at least one \emphd{full extension}, i.e., a full $k$-fold cover $\mathcal{H}' = (L, H')$ such that $H$ is a subgraph of $H'$.

    \begin{lem} \label{lem:full}
    Let $G$ be a $k$-critical graph for $k \geq 1$ and let $\mathcal{H} = (L,H)$ be a bad $(k-1)$-fold cover of $G$. If $\mathcal{H}$ is not full, then no full extension of $\mathcal{H}$ is canonical.
\end{lem}
\begin{scproof}
        Toward a contradiction, suppose $\mathcal{H}$ is not full and $\mathcal{H}_{F} = (L,H_{F})$ is a canonical full extension of $\mathcal{H}$. Then we can write $L(v) = \set{v_1,\ldots,v_{k-1}}$ for all $v \in V(G)$ so that \[E(H_F) \,=\, \set{u_iv_i \,:\, uv \in E(G),\, i \in [k-1]}.\] Since $\mathcal{H}$ is not full, there exist an edge $xy \in E(G)$ and an index $i \in [k-1]$ with $x_iy_i \notin E(H)$. As $G$ is $k$-critical, there is a proper $(k-1)$-coloring $c \colon V(G) \to [k-1]$ of $G - xy$. By permuting the colors if necessary, we may arrange that $c(x) = i$. This implies that $c(y) = i$ as well (otherwise $c$ would be a proper $(k-1)$-coloring of $G$), and hence \[T \,\defeq\, \set{u_{c(u)} \,:\, u \in V(G)}\] is a proper $\mathcal{H}$-coloring of $G$, a contradiction. 
\end{scproof}

    Lemma~\ref{lem:full} yields the following result, which shows that we may only consider full $(k-1)$-fold covers when trying to prove that a graph is robustly $k$-critical:

    \begin{cor}\label{corl:full}
        Let $G$ be a graph that is $k$-critical but not robustly $k$-critical for some $k \geq 1$. Then $G$ has a non-canonical bad full $(k-1)$-fold cover.
    \end{cor}
    \begin{scproof}
        Since $G$ is not robustly $k$-critical, it has a non-canonical bad $(k-1)$-fold cover $\mathcal{H}$. By Lemma~\ref{lem:full}, any full extension of $\mathcal{H}$ has the desired properties.
    \end{scproof}

    If $\mathcal{H} = (L,H)$ is a cover of a graph $G$ and $U \subseteq V(G)$, then we let $L(U) \defeq \bigcup_{u \in U} L(u)$. The next lemma is useful in inductive arguments for proving robust criticality.

    \begin{lem}\label{lemma:induction}
        Let $G$ be a graph and let $I \subseteq V(G)$ be an independent set. Suppose that for some $k \geq 2$, $\mathcal{H} = (L,H)$ is a bad full $(k-1)$-fold cover of $G$ and $\lambda \colon L(I) \to [k-1]$ is a function such that:
        \begin{itemize}
            \item for each $x \in I$, $\rest{\lambda}{L(x)}$ is a bijection from $L(x)$ to $[k-1]$, and
            \item if $c_1$, $c_2 \in L(I)$ have a common neighbor in $H$, then $\lambda(c_1) = \lambda(c_2)$.
        \end{itemize}
        If $G - I$ is robustly $(k-1)$-critical, then the cover $\mathcal{H}$ is canonical.
    \end{lem}
    \begin{scproof}
        We define a map $\lambda' \colon V(H) \to [k-1]$ as follows. If $c \in L(I)$, then we set $\lambda'(c) \defeq \lambda(c)$. If $c \in V(H) \setminus L(I)$ has a neighbor $c' \in L(I)$, then we let $\lambda'(c) \defeq \lambda(c')$ (note that this does not depend on the choice of $c'$ by the properties of the function $\lambda$). Finally, if $u \in V(G) \setminus I$ is a vertex with no neighbor in $I$, then we let $\rest{\lambda'}{L(u)}$ be an arbitrary bijection from $L(u)$ to $[k-1]$. Since $\mathcal{H}$ is full, this defines $\lambda'$ on all of $V(H)$. We claim that $\lambda'$ is a canonical labeling of $\mathcal{H}$. 

        It is easy to see that $\rest{\lambda'}{L(u)}$ is a bijection from $L(u)$ to $[k-1]$ for all $u \in V(G)$, and hence we may write $L(u) = \set{u_1, \ldots, u_{k-1}}$ so that $\lambda'(u_i) = i$. Suppose $\lambda'$ is not a canonical labeling. Since $\mathcal{H}$ is full, this means that there exist an edge $uv \in E(G)$ and indices $i \neq j$ such that $u_i v_j \in E(H)$. By the construction of $\lambda'$, it must be that $u$, $v \notin I$. Let $T \defeq \set{x_i \,:\, x \in I}$. For each $w \in V(G) \setminus I$, the only color in $L(w)$ that may have a neighbor in $T$ is $w_i$. Therefore, removing the neighborhood of $T$ from $H$ yields a cover $\mathcal{H'} = (L', H')$ of $G - I$ in which every list is of size at least $k-2$. If $T'$ were a proper $\mathcal{H}'$-coloring of $G - I$, then $T \cup T'$ would be a proper $\mathcal{H}$-coloring of $G$, a contradiction. Since $G - I$ is robustly $(k-1)$-critical, we conclude that $\mathcal{H}'$ is a canonical $(k-2)$-fold cover of $G- I$ by Lemma~\ref{lemma:excess}. Hence, $|L'(w)| = k - 2$ for all $w \in V(G) \setminus I$, which means that $w_i \notin L'(w)$ (note that this implies that every vertex in $V(G) \setminus I$ has a neighbor in $I$). In particular, $u_i \notin L'(u)$. But then the matching corresponding to the edge $uv$ in $\mathcal{H}'$ is not perfect (because $v_j$ has no neighbor in $H'$), which is impossible as $\mathcal{H}'$ is canonical.
\end{scproof} 

    A simple consequence of Lemma~\ref{lemma:induction} is Proposition~\ref{prop:examples}\ref{item:join}, which says that the join of a robustly critical graph and a clique is again robustly critical: 

    \begin{propcopy}{prop:examples}
        \begin{enumerate}
            \item[{\ref{item:join}}]  If $G$ is robustly critical, then so is $G \vee K_t$ for any $t \in \N$.
        \end{enumerate}
    \end{propcopy}
    \begin{scproof}
        It is enough to show that $G \vee K_1$ is robustly critical, since the general result would then follow by induction on $t$. Suppose $G$ is robustly $(k-1)$-critical and let $G' \cong G \vee K_1$ be formed by adding to $G$ a universal vertex $x$. Notice that $G' - x = G$. It is easy to see that $G'$ is $k$-critical. Now take any bad full $(k-1)$-fold cover $\mathcal{H} = (L,H)$ of $G'$ (we may assume $\mathcal{H}$ is full by Corollary~\ref{corl:full}). We need to argue that $\mathcal{H}$ is canonical. To this end, let $\lambda \colon L(x) \to [k-1]$ be an arbitrary bijection. Then the conditions of Lemma~\ref{lemma:induction} are satisfied for $I = \set{x}$, and hence $\mathcal{H}$ is canonical, as desired. 
    \end{scproof}


    The next corollary invokes Lemma~\ref{lemma:induction} with $|I| = 2$:

    \begin{cor}\label{corl:pair} 
        Let $G$ be a critical graph and let $x$, $y \in V(G)$ be two non-adjacent vertices with at most one common neighbor. If $G - x - y$ is robustly critical, then $G$ is robustly critical as well.
\end{cor}
\begin{scproof}
    Let $k \defeq \chi(G - x - y) +1$, so $G - x - y$ is robustly $(k-1)$-critical. Since $x$ and $y$ are not adjacent, we have $\chi(G) \leq k$, and since $G$ is critical, it follows that $\chi(G) = k$. Thus, we need to argue that $G$ is robustly $k$-critical. To this end, let $\mathcal{H} = (L,H)$ be a bad full $(k-1)$-fold cover of $G$ (we may assume $\mathcal{H}$ is full by Corollary \ref{corl:full}). We need to show $\mathcal{H}$ is canonical. Define a function $\lambda \colon L(x) \cup L(y) \to [k-1]$ by ensuring that $\rest{\lambda}{L(x)} \colon L(x) \to [k-1]$ and $\rest{\lambda}{L(y)} \colon L(y) \to [k-1]$ are bijections and, if $x$ and $y$ have a (unique) common neighbor $z \in V(G)$, then $\lambda(c_1) = \lambda(c_2)$ whenever $c_1 \in L(x)$ and $c_2 \in L(y)$ have a common neighbor in $L(z)$. Now we may apply Lemma \ref{lemma:induction} with $I = \set{x,y}$ to conclude that the cover $\mathcal{H}$ is canonical, as desired.
\end{scproof}   

    We are now ready to prove Proposition~\ref{prop:examples}\ref{item:Ekab}:

\begin{propcopy}{prop:examples}
        \begin{enumerate}
            \item[{\ref{item:Ekab}}]  The graphs $E_{k,a,b}$ from Example~\ref{exmp:E} are robustly critical.
        \end{enumerate}
    \end{propcopy}
\begin{scproof}
    Let $G \cong E_{k,a,b}$ and let $X_1$, $X_2$, $Y_1$, $Y_2$, $z$ be as in Example~\ref{exmp:E}. We will argue that $G$ is robustly $k$-critical by induction on $k$. In the base case $k = 3$, we have $G \cong E_{3,1,1} \cong C_5$, which is robustly $k$-critical by \cite[Lemma 12]{KMG21}. Now we suppose that $k > 3$.

    Let us first consider the case $b = k-2$. Then $|Y_2| = 1$, so we can let $z'$ be the unique vertex in $Y_2$. Without loss of generality, we may assume that $|X_1| \geq |X_2|$ (otherwise we may switch the roles of $z$ and $z'$ and of $X_1$ and $X_2$). Since $k > 3$, this implies that  $|X_1|$, $|Y_1| \geq 2$. Pick arbitrary vertices $x \in X_1$ and $y \in Y_1$ and notice that the graph $G - x - y \cong E_{k-1,a-1,k-3}$ is robustly $(k-1)$-critical by the inductive hypothesis. Since $x$ and $y$ have a unique common neighbor, namely $z$, it follows that $G$ is robustly $k$-critical by Corollary~\ref{corl:pair}.

    Since $E_{k,a,b} \cong E_{k,b,a}$, we may now assume that $a$, $b \leq k - 3$. Since $a + b \geq k - 1$ by the definition of $E_{k,a,b}$, we have $a$, $b \geq 2$. If $a + b \geq k$, then we can pick arbitrary vertices $x \in X_1$, $y \in Y_1$, notice that $G - x - y \cong E_{k-1, a-1, b-1}$ is robustly $(k-1)$-critical by the inductive hypothesis, and apply Corollary~\ref{corl:pair} to conclude that $G$ is robustly $k$-critical, as desired. Therefore, we may assume that $a + b = k-1$, and hence $|X_1| = |Y_2| = a$ and $|X_2| = |Y_1| = b$.
    
    Let $\mathcal{H} = (L,H)$ be a bad full $(k-1)$-fold cover of $G$. 
    We need to argue that $\mathcal{H}$ is canonical. Let $L(z) = \set{z_1, \ldots, z_{k-1}}$ and pick arbitrary vertices $x \in X_1$ and $y \in Y_1$. Since $\mathcal{H}$ is full, we can write $L(x) = \set{x_1, \ldots, x_{k-1}}$ and $L(y) = \set{y_1, \ldots, y_{k-1}}$ so that $x_i z_i$, $y_i z_i \in E(H)$ for all $i \in [k-1]$. 

    Fix an arbitrary vertex $w \in Y_2$ and let $L(w) = \set{w_1, \ldots, w_{k-1}}$ where $w_i y_i \in E(H)$ for all $i \in [k-1]$. Define $\lambda \colon L(x) \cup L(w) \to [k-1]$ by $\lambda(x_i) = \lambda(w_i) \defeq i$. Note that the graph $G - x - w \cong E_{k-1, a-1, b}$ is robustly $(k-1)$-critical by the inductive hypothesis. Hence, if we could apply Lemma~\ref{lemma:induction} to the independent set $I = \set{x, w}$ and the function $\lambda$, it would follow that $\mathcal{H}$ is canonical, as desired. On the other hand, if Lemma~\ref{lemma:induction} cannot be applied, then there exist colors $x_i \in L(x)$ and $w_j \in L(w)$ with $i \neq j$ that have a common neighbor $c \in L(v)$ for some $v \in X_2$. We claim that in this case $G$ is $\mathcal{H}$-colorable, which is a contradiction. Indeed, we can construct a proper $\mathcal{H}$-coloring of $G$ as follows. From $L(x)$ and $L(y)$, we pick the colors $x_i$ and $y_i$ respectively. Removing the neighbors of $x_i$ and $y_i$ from $L(u)$ for all $u \in X_2 \cup Y_2$ yields a $(k-2)$-fold cover $\mathcal{H}'$ of the $(k-1)$-clique $G[X_2 \cup Y_2]$. This cover is not canonical because the matching corresponding to the edge $vw$ is not prefect (the color $w_j$ has no neighbor in it). Since cliques are robustly critical (for example, by Proposition~\ref{prop:examples}\ref{item:join}), there exists a proper $\mathcal{H}'$-coloring $T'$ of $G[X_2 \cup Y_2]$. We can then extend $T' \cup \set{x_i, y_i}$ to an independent transversal $T$ of $\mathcal{H}$ by picking colors for the remaining vertices of $G$ greedily in the following order: the vertices in  $X_1 \setminus \set{x}$, then the ones in $Y_1 \setminus \set{y}$, and finally $z$.
\end{scproof}

    \section{Joins with large cliques are strongly chromatic-choosable: Proof of Theorem~\ref{theo:vertex}}\label{sec:vertex}

    In this section we show that the join of a vertex-critical graph and a large enough clique is strongly chromatic-choosable. The starting point of our analysis is the following fact established by Kostochka, Zhu, and the first named author:

    \begin{thm}[{AB--Kostochka--Zhu \cite[Theorem 1.5]{BK17}}]\label{theo:DP-chrom-choosable1}
        Let $G$ be a graph with $m$ edges and let $k \defeq \chi(G)$. Then $\chi_{\mathsf{DP}}(G \vee K_t) = \chi(G \vee K_t) = k + t$ for all $t \geq 3m$.
    \end{thm}

    We also need the following slightly more technical variant of Theorem~\ref{theo:DP-chrom-choosable1} that is well-suited for inductive arguments: 

    \begin{thm}[{AB--Kostochka--Zhu \cite[Theorem 2.1]{BK17}}]\label{theo:DP-chrom-choosable}
        Let $G$ be a graph with $\chi(G) \leq k$ and let $J \cong G \vee K_t$ be the graph obtained from $G$ be adding a $t$-element set $U$ of universal vertices. Suppose that $\mathcal{H} = (L,H)$ is a cover of $J$ such that $|L(u)| \geq k + t$ for each $u \in U$ and 
        \[
            t \,\geq\, \frac{3}{2} \sum_{v \in V(G)} \max \big\{\deg_G(v) + t - |L(v)| + 1, \, 0 \big\}.
        \]
        Then $J$ is $\mathcal{H}$-colorable.
    \end{thm}
    
    We use Theorem~\ref{theo:DP-chrom-choosable} and induction on $t$ to prove the following statement that implies Theorem~\ref{theo:vertex}:

    \begin{lem} \label{lem: bigcomplete}
        Let $G$ be a $k$-vertex-critical graph and let $J \cong G \vee K_t$ be the graph obtained from $G$ be adding a $t$-element set $U$ of universal vertices. Suppose that $L$ is a list assignment for $J$ such that $|L(x)| \geq k + t - 1$ for each $x \in V(J)$ and 
        \begin{equation}\label{eq:t_bound}
            t \,\geq\, \frac{3}{2} \sum_{v \in V(G)} \max \big\{\deg_G(v) + t - |L(v)| + 1, \, 0 \big\}.
        \end{equation}
        If $L$ is not a constant $(k + t - 1)$-assignment, then $J$ is $L$-colorable.
    \end{lem}
\begin{scproof}
    The proof is by induction on $t$.  When $t=0$, \eqref{eq:t_bound} yields that $|L(v)| \geq \deg_G(v) + 1$ for each $v \in V(G)$, which implies that $J = G$ is $L$-colorable (using a greedy coloring procedure).

    Now suppose $t > 0$ and the statement is true for all smaller values of $t$. Since $G$ is $k$-vertex-critical, for any $v \in V(G)$, $\chi(G - v) = k - 1$ and hence the graph $J - v$ is $L$-colorable by Theorem~\ref{theo:DP-chrom-choosable}. 

    We may assume that $\deg_G(v) + t - |L(v)| + 1 > 0$ for all $v \in V(G)$. Indeed, suppose there is a vertex $v \in V(G)$ such that $\deg_G(v) + t + 1 \leq |L(v)|$.  As discussed above, there is a proper $L$-coloring of $J - v$.  Since $v$ is adjacent to $\deg_G(v) + t < |L(v)|$ vertices in $J$, this coloring can be extended to a proper $L$-coloring of $J$, as desired. 

        Suppose that for some $u_0 \in U$ there exist a color $c_0 \in L(u_0)$ and a vertex $v_0 \in V(G)$ such that $c_0 \notin L(v_0)$. 
        Let $U' \defeq U \setminus \{u_0\}$ and $J' \defeq J - u_0 \cong G \vee K_{t-1}$, and let $L'$ be the list assignment for $J'$ given by $L'(x) \defeq L(x) \setminus \set{c_0}$ for each $x \in V(J')$.  Notice that for all $v \in V(G) \setminus \{v_0\}$,
        \[
            \deg_G(v) + (t-1) - |L'(v)| + 1 \,\leq\, \deg_G(v) + t - |L(v)| + 1.
        \]
        Furthermore, $\deg_G(v_0) + (t-1) - |L'(v_0)| + 1 = \deg_G(v_0) + t - |L(v_0)|$ as $L'(v_0) = L(v_0)$.  Thus, 
        \[
            t - 1 \,\geq\, \frac{3}{2} \sum_{v \in V(G)} \max \big\{\deg_G(v) + (t-1) - |L'(v)| + 1, \,0 \big\}.
        \]
        We also have $|L'(x)| \geq k + t -2$ for all $x \in V(J')$. Moreover, $|L'(v_0)| = |L(v_0)| \geq k + t - 1$, so $L'$ is not a constant $(k + t -2)$-assignment. Therefore, by the inductive hypothesis, $J'$ is $L'$-colorable. Any proper $L'$-coloring of $J'$ can be extended to a proper $L$-coloring of $J$ by giving the color $c_0$ to $u_0$. It follows that $J$ is $L$-colorable, as desired. Thus, letting $\mathcal{L} \defeq \bigcup_{u \in U} L(u)$, we may assume that $\mathcal{L} \subseteq L(v)$ for all $v \in V(G)$. 

        Note that $|\mathcal{L}| \geq t + k - 1$. If $|\mathcal{L}| = t + k - 1$, then $L(u) = \mathcal{L}$ for all $u \in U$. Since $L$ is not a constant $(k + t - 1)$-assignment, it follows that there exists a vertex $v_0 \in V(G)$ such that $\mathcal{L}$ is a proper subset of $L(v_0)$. Consider the graph $J_0 \defeq J - v_0 \cong (G - v_0) \vee K_t$. Since $G$ is $k$-vertex-critical, $\chi(J_0) = k + t - 1$, so we can fix a proper $L$-coloring of $J_0$ using only the colors from $\mathcal{L}$. We can then extend it to a proper $L$-coloring of $G$ by assigning to $v_0$ an arbitrary color from $L(v_0) \setminus \mathcal{L}$.
        
        Thus, we may assume $|\mathcal{L}| \geq t + k$. We can now form a proper $L$-coloring of $J$ in two stages.  First, we greedily color the vertices in $U$ using $t$ colors (this is possible since $|L(u)| \geq k + t-1 \geq t$ for each $u \in U$). Let $\mathcal{C} \subseteq \mathcal{L}$ be the set of colors used on $U$. Notice that $|\mathcal{L} \setminus \mathcal{C}| = |\mathcal{L}| - t \geq k$ and $\mathcal{L} \setminus \mathcal{C} \subseteq L(v)$ for each $v \in V(G)$. Since $G$ is $k$-colorable, we can now properly color $G$ using the colors from $\mathcal{L} \setminus \mathcal{C}$, which completes a proper $L$-coloring of $J$ and finishes the proof. 
\end{scproof}


    \begin{theocopy}{theo:vertex}
        If $G$ is a vertex-critical graph with $m$ edges, then for all $t \geq 3m$, the graph $G \vee K_t$ is strongly chromatic-choosable.
    \end{theocopy}
\begin{scproof}
    Let $k \defeq \chi(G)$. We may assume $k \geq 2$ (since the only $1$-vertex-critical graph is $K_1$ and cliques are strongly chromatic-choosable). The graph $G \vee K_t$ is easily seen to be $(k + t)$-vertex-critical. Now let $L$ be a bad $(k + t - 1)$-assignment for $G \vee K_t$. Then 
    \begin{align*}
        \frac{3}{2} \sum_{v \in V(G)} \max \big\{\deg_G(v) + t - |L(v)| + 1, \, 0 \big\} \,&=\, \frac{3}{2} \sum_{v \in V(G)} \max \big\{\deg_G(v) - k + 2,\, 0\big\} \\
        &\leq\, \frac{3}{2} \sum_{v \in V(G)} \deg_G(v) \,=\, 3m.
    \end{align*}
    Therefore, if $t \geq 3m$, we may apply Lemma~\ref{lem: bigcomplete} to conclude that $G \vee K_t$ is $L$-colorable unless $L$ is constant, as desired.
\end{scproof}

    \section{Joins with larger cliques are robustly critical: Proof of Theorem~\ref{theo:joins_make_robust}}\label{sec:robust}

    The bulk of our proof of Theorem~\ref{theo:joins_make_robust} is in the following lemma, in which we assume that the graph $G$ satisfies $\chi_\mathsf{DP}(G) = \chi(G)$:

    \begin{lem} \label{lem: reduce}
        If $G$ is an $n$-vertex $k$-critical graph with $\chi_{\mathsf{DP}}(G)=k$, then $G \vee K_t$ is robustly critical for all $t \geq n^3$. 
    \end{lem}
    \begin{scproof}
        We may assume that $k \geq 2$ (the only $1$-critical graph is $K_1$, and cliques are robustly critical by Proposition~\ref{prop:examples}\ref{item:join}), and hence $n \geq 2$. Let $J \cong G \vee K_t$ be the graph obtained from $G$ be adding a $t$-element set $U$ of universal vertices, where $t \geq n^3$. It is easy to see that $J$ is $(k + t)$-critical. Let $\mathcal{H} = (L,H)$ be a bad full $(k + t - 1)$-fold cover of $J$. Our goal is to show that $\mathcal{H}$ is canonical.

        We begin with some terminology and notation. An \emphd{independent partial transversal} in $\mathcal{H}$ is an independent set $T \subseteq V(H)$ such that $|L(v) \cap T| \leq 1$ for all $v \in V(J)$. 
        Given an independent partial transversal $T$, we define the \emphd{domain} of $T$ as
        $
            \dom(T) \defeq \set{v \in V(J) \,:\, L(v) \cap T \neq \0}
        $
        and let $\mathcal{H}_T = (L_T, H_T)$ be the cover of the graph $J - \dom(T)$ obtained by removing the neighborhood of $T$ from $H$. 
        Note that if $T'$ is a proper $\mathcal{H}_T$-coloring of $J - \dom(T)$, then $T \cup T'$ is a proper $\mathcal{H}$-coloring of $J$. It follows that $\mathcal{H}_T$ must be a bad cover of $J - \dom(T)$. 
        
        Let $\mathcal{T}_U$ be the set of all independent partial transversals $T$ with $\dom(T) \subseteq U$. The \emphd{excess} of $T \in \mathcal{T}_U$ is the quantity $\epsilon(T)$ defined as follows:
        \[
            \epsilon(T) \,\defeq\, \sum_{v \in V(G)} \big(|L_T(v)|-(k + t - 1 - |T|)\big).
        \]
        All the above definitions make sense when $T = \0$. In particular, $\mathcal{H}_\0 = \mathcal{H}$ and $\epsilon(\0) = 0$.

        Notice that if $T \subseteq T' \in \mathcal{T}_U$, then $\epsilon(T) \leq \epsilon(T')$. We define a sequence $T_0 \subset T_1 \subset \cdots \subset T_N \in \mathcal{T}_U$ of independent partial transversals with strictly increasing excess. Start by setting $T_0 \defeq \0$. Once $T_i$ is constructed, we proceed as follows.

        \begin{leftbar}
            \noindent \textbf{If} there exist vertices $u$, $u' \in U \setminus \dom(T_i)$, $v \in V(G)$ and colors $c \in L_{T_i}(u)$, $c' \in L_{T_i}(u')$ such that 
            \[
                u \neq u', \qquad cc' \notin E(H), \qquad \text{and} \qquad N_H(c) \cap N_H(c') \cap L_{T_i}(v) \neq \0,
            \]
            then we pick any such $u$, $u'$, $v$, $c$, $c'$ and let $T_{i + 1} \defeq T \cup \set{c, c'}$. 

            \medskip

             \noindent \textbf{Else, if} there exist vertices $u \in U \setminus \dom(T_i)$, $v \in V(G)$ and a color $c \in L_{T_i}(u)$ with
             \[
                N_H(c) \cap L_{T_i}(v) \,=\, \0,
            \]
            then we pick any such $u$, $v$, $c$ and let $T_{i+1} \defeq T_i \cup \set{c}$. 

            \medskip

             \noindent \textbf{Otherwise}, we let $N \defeq i$ and terminate the construction.
        \end{leftbar}

        \noindent Let $T \defeq T_N$ and define
        \begin{align*}
            D \,&\defeq\, \dom(T), \qquad U' \,\defeq\, U \setminus D, \qquad t' \,\defeq\, |U'| \,=\, t - |T|,\\
            J' \,&\defeq\, J - D \,\cong\, G \vee K_{t'}, \qquad \text{and} \qquad \mathcal{H}' \,=\, (L', H') \,\defeq\, \mathcal{H}_T.
        \end{align*}
        Note that for all $x \in V(J')$, \[|L'(x)| \,\geq\, |L(x)| - |T| \,=\,k +t' - 1.\]
        Let us record some basic properties of the above construction:

        \begin{cl} \label{cl: algorithm}
        The following statements are valid.

        \begin{enumerate}[label=\ep{\normalfont\roman*}]
            \item\label{item:size_excess} For all $0 \leq i \leq N$, we have $|T_i| \leq 2i$ and $\epsilon(T_i) \geq i$.

            \item\label{item:common_nbr} If $u$, $u' \in U'$ are distinct vertices and $c \in L'(u)$, $c' \in L'(u')$ are colors that have a common neighbor in $L'(V(G))$, then $cc' \in E(H')$.

            \item\label{item:has_nbr} For all $u \in U'$ and $v \in V(G)$, every color in $L'(u)$ has a neighbor in $L'(v)$.
        \end{enumerate}
    \end{cl}
    \begin{claimproof}
        Statement \ref{item:size_excess} follows since $|T_{i+1}| \leq |T_i| + 2$ and $\epsilon(T_{i+1}) \geq \epsilon(T_i) + 1$ for all $i < N$. Statements \ref{item:common_nbr} and \ref{item:has_nbr} hold because the process terminated on step $N$.
\end{claimproof}

    We will eventually show that $N = 0$ and $T = \0$, i.e., the above process terminates immediately (this must be the case if $\mathcal{H}$ is canonical which is what we are aiming for). Before we can do that, we need to establish a few claims about the structure of the cover $\mathcal{H}'$. We start with an upper bound on the excess of $T$:

\begin{cl} \label{cl: excess}
    $\epsilon(T) \leq n(n-k)$.
\end{cl}    
\begin{claimproof}
    Suppose not and let $0 < i \leq N$ be the minimum index such that $\epsilon(T_i) > n(n-k)$. Then $\epsilon(T_{i-1}) \leq n (n-k)$, so, by Claim~\ref{cl: algorithm}\ref{item:size_excess}, $|T_{i}| \leq |T_{i - 1}| + 2 \leq 2n (n-k) + 2$. Set $t_i \defeq t - |T_i|$. Then 
    \begin{equation}\label{eq:t'}
        t_i \,\geq\, n^3 - 2n(n-k) - 2 \,\geq\, n^3 - 2n(n-2) - 2 \,\geq\,  3n(n-1)/2,
    \end{equation}
    where the last inequality is valid for all $n \geq 1$. 
    As $\epsilon(T_i) > n(n-k)$, there must be a vertex $z \in V(G)$ with $|L_{T_i}(z)| - (k+t - 1-|T_i|) > n - k$, i.e., $|L_{T_i}(z)| \geq n + t_i$. Consider the graph \[J - \dom(T_i) - z \,\cong\, (G - z) \vee K_{t_i}.\] Since $G$ is $k$-critical, we have $\chi(G-z) = k-1$. Inequality \eqref{eq:t'} implies that $t_i \geq 3|E(G-z)|$, so, by Theorem~\ref{theo:DP-chrom-choosable1}, $\chi_{\mathsf{DP}}(J  - \dom(T_i) - z) = k+t_i - 1$.  Since $|L_{T_i}(v)| \geq |L(v)| - |T_i| = k+t_i - 1$ for each $v \in V(J - \dom(T_i) - z)$, there is a proper $\mathcal{H}_{T_i}$-coloring $T'$ of $J-z - \dom(T_i)$.  As $z$ has at most $n-1+t_i < |L_{T_i}(z)|$ neighbors in $J - \dom(T_i)$, we can find an element in $L_{T_i}(z)$ to add to $T \cup T'$ in order to complete a proper $\mathcal{H}$-coloring of $J$. This contradicts the fact that $\mathcal{H}$ is a bad cover.
\end{claimproof}

    Now we show that not every vertex of $G$ can be making a positive contribution to the excess of $T$: 

    \begin{cl} \label{cl: noexcess}
        There is a vertex $v \in V(G)$ satisfying $|L'(v)|-(k+t' - 1)=0$.
    \end{cl}
    \begin{claimproof}
        This is the only place in the proof where we use the assumption that $\chi_\mathsf{DP}(G) = k$. Suppose $|L'(v)| \geq k+t'$ for all $v \in V(G)$.  Since $|L'(u)| \geq k+t' - 1 > t'$ for each $u \in U'$, we can greedily form a proper $\mathcal{H}'$-coloring $T'$ of the clique $J[U']$. As $|L_{T \cup T'}(v)| \geq |L'(v)| - t' \geq k$ for each $v \in V(G)$ and $\chi_{\mathsf{DP}}(G) = k$, there is a proper $\mathcal{H}_{T \cup T'}$-coloring of $G$, and thus $J$ is $\mathcal{H}$-colorable, a contradiction.
    \end{claimproof}

    We can now show that the cover $\mathcal{H}'$ is canonical when restricted to $U'$:

    \begin{cl} \label{cl: canonical}
    For all $u \in U'$, $|L'(u)| = k + t' - 1$. Furthermore, it is possible to list the colors in $L'(u)$ as $L'(u) = \set{u_1, \ldots, u_{k+t'-1}}$ so that $u_i u_i' \in E(H)$ for all distinct $u$, $u' \in U'$ and all $i \in [k + t' - 1]$. 
\end{cl}  
    \begin{claimproof}
        Let $v \in V(G)$ be a vertex with $|L'(v)| = k + t' - 1$ given by Claim~\ref{cl: noexcess}.  For any $u \in U'$, on the one hand, we have $|L'(u)| \geq k+t'-1$. On the other hand, by Claim~\ref{cl: algorithm}\ref{item:has_nbr}, each color in $L'(u)$ has a neighbor in $L'(v)$, so $|L'(u)| \leq |L'(v)| = k + t' - 1$. It follows that $|L'(u)| = k + t' - 1$ and the matching $E_{H'}(L'(u), L'(v))$ is perfect. Thus, we can write $L'(v) = \set{v_1, \ldots, v_{k + t' - 1}}$ and $L'(u) = \set{u_1, \ldots, u_{k + t' - 1}}$ so that $u_i v_i \in E(H)$ for all $i \in [k+t' - 1]$.  By Claim~\ref{cl: algorithm}\ref{item:common_nbr}, if $u$, $u' \in U'$ are distinct vertices, then $u_i u_i' \in E(H)$ for all $i \in [k + t' - 1]$, as claimed.
    \end{claimproof}



    For each $c \in L'(V(G))$, let $M(c)$ be the set of all $u \in U'$ such that $c$ has no neighbor in $L'(u)$.

\begin{cl} \label{cl: nomatch}
For each  $v \in V(G)$ and $c \in L'(v)$, $|M(c)| < 3(n-k+1)(n-1)/2 \leq 3(n-1)^2/2$.
\end{cl}
\begin{claimproof}
    Toward a contradiction, suppose $v \in V(G)$ and $c \in L'(v)$ satisfy $|M(c)| \geq 3(n-k+1)(n-1)/2$. Set $U^* \defeq U' \setminus M(c)$. Since $|L'(u)| = k + t' - 1 > t' > |U^*|$ for all $u \in U'$, we can greedily form a proper $\mathcal{H}'$-coloring $T'$ of the clique $J[U^*]$ so that no color in $T'$ is adjacent to $c$. We will argue that $T \cup T' \cup \set{c}$ can be extended to a proper $\mathcal{H}$-coloring of $J$, contradicting the fact that $\mathcal{H}$ is bad.

    To this end, consider the graph \[J'' \,\defeq\, J' - U^* - v \,\cong\, (G - v) \vee K_{|M(c)|},\] and let $\mathcal{H}'' = (L'', H'')$ be the cover of $J''$ obtained by removing from $H'$ the neighborhood of $T' \cup \set{c}$. We show that $J''$ is $\mathcal{H}''$-colorable. By the definition of $M(c)$, every vertex $u \in M(c)$ satisfies
    \[
        |L''(u)| \,\geq\, |L'(u)| - |T'| \, =\, k + |M(c)| - 1.
    \]
    Moreover, for each $x \in V(G) \setminus \set{v}$,
    \[
        |L''(x)| \,\geq\, |L'(x)| - |T'| - 1 \,\geq\, k + |M(c)| - 2.
    \]
    It follows that
    \begin{align*}
        &\frac{3}{2} \sum_{x \in V(G) \setminus \set{v}} \max \big\{\deg_{G - v} (x) + |M(c)| - |L''(x)| + 1, \, 0\big\} \\
        \leq\,& \frac{3}{2} \sum_{x \in V(G) \setminus \set{v}}  \max \big\{(n-2) + |M(c)| - (k + |M(c)| - 2) + 1, \, 0 \big\} \\
        =\,& \frac{3(n-k+1)(n-1)}{2} \,\leq\, |M(c)|.
    \end{align*}
    Since $G$ is $k$-critical, $\chi(G - v) = k - 1$, so we may apply Theorem~\ref{theo:DP-chrom-choosable} with $G - v$ in place of $G$, $k-1$ in place of $k$, and $|M(c)|$ in place of $t$ to conclude that $J''$ is $\mathcal{H}''$-colorable, as desired.
\end{claimproof}

 We now use the above claims to conclude that $T = \0$:

\begin{cl} \label{cl: zeroexcess}
    $T = \0$.
\end{cl}
\begin{claimproof}
    By Claim~\ref{cl: algorithm}\ref{item:size_excess}, it suffices to argue that $\epsilon(T) = 0$, i.e., $|L'(v)| = k+t' - 1$ for all $v \in V(G)$. Suppose for contradiction that there is a vertex $z \in V(G)$ such that $|L'(z)| \geq k+t'$. It follows from Claims~\ref{cl: algorithm}\ref{item:size_excess} and \ref{cl: excess} that $|T| \leq 2 n (n-k) \leq 2n(n-2)$, and hence $t' = t - |T| \geq n^3 - 2n(n-2)$. Therefore, by Claim~\ref{cl: nomatch}, for each $c \in L'(z)$, we have
    \begin{align*}
        \frac{t'}{2} - |M(c)| \,\geq\, \frac{n^3 - 2n(n-2)}{2} - \frac{3(n-1)^2}{2} \,>\, 0,
    \end{align*}
    where the last inequality holds for all $n \geq 1$. 

    Using Claim~\ref{cl: canonical}, we write $L'(u) = \set{u_1, \ldots, u_{k + t' - 1}}$ for each $u \in U'$ so that $u_i u_i' \in E(H)$ for all distinct $u$, $u' \in U'$ and $i \in [k + t' - 1]$.  Define a function $f \colon L'(z) \to [k + t' - 1]$ by making $f(c) = i$ if $cu_i \in E(H)$ for some $u \in U' \setminus M(c)$. By Claim~\ref{cl: algorithm}\ref{item:common_nbr}, the value $f(c)$ does not depend on the choice of the vertex $u \in U' \setminus M(c)$. Since $|L'(z)| \geq k + t'$, by the Pigeonhole Principle, there exist distinct colors $c_1$, $c_2 \in L'(z)$ and an index $i \in [k + t' - 1]$ such that $f(c_1) = f(c_2) = i$. Since $|M(c_1)| + |M(c_2)| <  t'$, there is a vertex $u \in U' \setminus (M(c_1) \cup M(c_2))$. Then, by the definition of $f$, both $c_1$ and $c_2$ are adjacent to $u_i$, which is impossible as $E_H(L(z), L(u))$ is a matching.
\end{claimproof}

   If follows from Claim~\ref{cl: zeroexcess} that $U' = U$, $t' = t$, $J' = J$, and $\mathcal{H}' = \mathcal{H}$. We are now ready to perform the final step in the proof of Lemma~\ref{lem: reduce}.
   
   Using Claim~\ref{cl: canonical}, we write $L(u) = \set{u_1, \ldots, u_{k + t -1}}$ for each $u \in U$ so that $u_i u_i' \in E(H)$ for all distinct $u$, $u' \in U$ and $i \in [k + t - 1]$. For every $v \in V(G)$, we let $L(v) = \set{v_1, \ldots, v_{k+t - 1}}$ so that for all $u \in U$ and $i \in [k + t - 1]$, $u_iv_i \in E(H)$; this is possible because the cover $\mathcal{H}$ is full and thanks to Claim~\ref{cl: algorithm}\ref{item:common_nbr}. We will show that the mapping $\lambda \colon V(H) \to [k + t - 1]$ given by $\lambda(x_i) \defeq i$ for all $x \in V(J)$ is a desired canonical labeling of $\mathcal{H}$. To this end, we only need to argue that $x_iy_i \in E(H)$ for all $i \in [k + t - 1]$ and $xy \in E(G)$, since for the other edges of $J$ this holds by construction.

    Let $\Gamma$ be the simple graph with vertex set $[k+t - 1]$ in which distinct $i$, $j \in [k+t - 1]$ are adjacent if and only if there is an edge $uv \in E(G)$ such that $u_i v_j \in E(H)$ or $u_j v_i \in E(H)$.  Notice that each edge of $G$ contributes at most $2$ to the degree of each $i \in V(\Gamma)$.  Consequently, $\Delta(\Gamma) \leq 2|E(G)|$. A greedy construction starting with an arbitrary vertex shows that every vertex of $\Gamma$ belongs to an independent set of size at least 
    \[
        \frac{|V(\Gamma)|}{1 + \Delta(\Gamma)} \,\geq\,\frac{k+t - 1}{1 + 2|E(G)|} \,\geq\, \frac{k+t - 1}{1 + n(n-1)}.
    \]
Since $t \geq n^3 > n(n-1)(k-1)$, the last quantity is at least $k-1$. 

Consider any independent set $I \subseteq [k + t - 1]$ of size $k-1$ in $\Gamma$. Fix an arbitrary proper coloring $c$ of the clique $J[U]$ using the colors in the $t$-element set $[k+t - 1] \setminus I$ and let
    \[P \,\defeq\, \{u_{c(u)} \,:\, u \in U \}.\]
Then $P$ is an independent partial transversal of $\mathcal{H}$, and the cover $\mathcal{H}_P = (L_P, H_P)$ of $G$ is $(k-1)$-fold and satisfies $L_P(v) = \set{v_i \,:\, i \in I}$ for all $v \in V(G)$. Since $I$ is independent in $\Gamma$, for all $xy \in E(G)$, 
    \begin{equation}\label{eq:almost_canonical}
        E_H(L_P(x), L_P(y)) \,\subseteq\, \set{x_i y_i \,:\, i \in I}.
    \end{equation}
    We claim that the cover $\mathcal{H}_P$ is full, i.e., the inclusion in \eqref{eq:almost_canonical} is actually an equality. Otherwise, there exist an edge $xy \in E(G)$ and an index $i \in I$ with $x_i y_i \notin E(H)$. Since $G$ is $k$-critical, we can fix a proper $(k-1)$-coloring $c'$ of $G - xy$ using the color set $I$. Furthermore, by permuting the colors, we may assume $c'(x) = i$, which implies that $c'(y) = i$ as well. Since $x_i y_i \notin E(H)$, it follows that \[P' \,\defeq\, \set{v_{c'(v)} \,:\, v \in V(G)}\] is a proper $\mathcal{H}_P$-coloring of $G$, and hence $P \cup P'$ is a proper $\mathcal{H}$-coloring of $J$, a contradiction.

    To summarize, for all $i \in I$ and $xy \in E(G)$, we have $x_iy_i \in E(H)$. Since every $i \in [k + t - 1]$ belongs to some $(k-1)$-element independent set $I$ in $\Gamma$, we conclude that the labeling $\lambda$ given by $\lambda(x_i) \defeq i$ for all $x \in V(J)$ and $i \in [k + t - 1]$ is indeed canonical, and the proof is complete.
\end{scproof}


\begin{theocopy}{theo:joins_make_robust}
        If $G$ is a critical graph with $m$ edges, then for all $t \geq 100\,m^3$, the graph $G \vee K_t$ is robustly critical.
    \end{theocopy}
 \begin{scproof}
    Let $k \defeq \chi(G)$. We may assume $k \geq 3$, since otherwise $G$ is a complete graph and we are done by Proposition~\ref{prop:examples}\ref{item:join}. In particular, we have $|V(G)| \leq m$. Note that \[G \vee K_t \,\cong\, G \vee K_{3m} \vee K_{t - 3m}.\] Since $G$ is $k$-critical, $G \vee K_{3m}$ is $(k+3m)$-critical and, by Theorem~\ref{theo:DP-chrom-choosable1}, $\chi_{\mathsf{DP}}(G \vee K_{3m}) = k +3m$. Thus, by Lemma~\ref{lem: reduce}, $G \vee K_t$ is robustly critical provided that $t - 3m \geq |V(G \vee K_{3m})|^3$. This is indeed the case as
    \[
        |V(G \vee K_{3m})|^3 \,\leq\, (|V(G)| + 3m)^3 \,\leq\, 64 m^3 \,\leq\, 100 m^3 - 3m \,\leq\, t - 3m. \qedhere
    \]
\end{scproof}

\subsection*{Acknowledgments}

We are grateful to the anonymous referees for their helpful feedback. This material is based upon work partially supported by the Alfred P. Sloan Foundation and the National Science Foundation under grant DMS-2528522. Any opinions, findings, and conclusions or recommendations expressed in this material are those of the authors and do not necessarily reflect the views of the National Science Foundation.

    \printbibliography

\end{document}